\documentclass[jimaging,article,submit,moreauthors,pdftex]{Definitions/mdpi}

\firstpage{1} 
\pubvolume{xx}
\issuenum{1}
\articlenumber{5}
\pubyear{2020}
\copyrightyear{2020}
\history{Received: date; Accepted: date; Published: date}

\pdfoutput=1

\preto{\abstractkeywords}{\nolinenumbers}

\usepackage{algorithmic}
\usepackage{algorithm}
\usepackage{subfig}

\Title{The Empirical Watershed Wavelet}

\Author{Basile Hurat $^{1,2}$, Zariluz Alvarado$^{1,3}$ and J\'er\^ome Gilles$^{1,4,}$*\orcidC{}}

\AuthorNames{Basile Hurat, Zariluz Alvarado and J\'er\^ome Gilles}

\address{%
$^{1}$ \quad Department of Mathematics \& Statistics, San Diego State University, 5500 Campanile Dr, San Diego, CA 92182-7720\\
$^{2}$ \quad bhurat2855@sdsu.edu\\
$^{3}$ \quad zalvarado3702@sdsu.edu\\
$^{4}$ \quad jgilles@sdsu.edu
}
\corres{Correspondence: jgilles@sdsu.edu; Tel.: +1-619-594-7240}

\abstract{The empirical wavelet transform is an adaptive multiresolution analysis tool based on the idea of building filters on a data-driven partition of the Fourier domain. However, existing 2D extensions are constrained by the shape of the detected partitioning. In this paper, we provide theoretical results that permits us to build 2D empirical wavelet filters based on 
an arbitrary partitioning of the frequency domain. We also propose an algorithm to detect such partitioning from an image spectrum by combining a scale-space representation to estimate the position of dominant harmonic modes and a watershed transform to find the boundaries of the different supports making the expected partition. This whole process allows us to define the 
empirical watershed wavelet transform. We illustrate the effectiveness and the advantages of such adaptive transform, first visually on toy images, and next on both unsupervised texture segmentation and image deconvolution applications.}


\begin{document}


\section{Introduction}
\label{sec:1}
Within the field of image processing and computer vision, multi-resolution analysis is a vastly applicable tool. For example, the use of multi-resolution analysis and compressive sensing theory has led to state of the art denoising and deconvolution techniques. Also, the use of multi-resolution analysis has shown to be quite effective in texture analysis. 

Wavelets have been the standard tool for multi-resolution analysis, and their effectiveness is matched by their popularity. However, wavelets are prescriptive in their construction in the sense that the wavelet filters are built independently of the data. Adaptive (i.e. data-driven) methods have shown many promising improvements in applications across all fields of science. Hereafter, we review the most popular adaptive methods developed in the last two decades.

An early adaptive method is the empirical mode decomposition (EMD) \cite{emd1}, \cite{emd2}. The EMD is an iterative algorithmic method that attempts to decompose a signal into its amplitude-modulated frequency-modulated components (also called harmonic modes). Despite its success in analyzing non-stationary signals, the EMD lacks solid mathematical foundation due to its algorithmic nature, making difficult to predict its behavior. Furthermore, the bidimensional empirical mode decomposition is computationally inefficient, especially with the presence of noise \cite{emdproblem}. 
An alternative, called variational mode decomposition (VMD), was proposed by Dragomiretskiy et al. \cite{vmd1}, \cite{vmd2}. It aims at decomposing the spectrum of an image into a mixture of Gaussians with specific constraints, using optimization techniques. While this method has seen success in many applications \cite{vmdapp4}, \cite{vmdapp2}, \cite{vmdapp5}, \cite{vmdapp1}, \cite{vmdapp3}. the assumption that a spectrum can be meaningfully represented as a mixture of Gaussians can lead to a loss of information. 
Yet another approach is using synchrosqueezed wave packets, proposed by Daubechies et al. in \cite{synchrosqueeze1}, extended into 2D by Yang et al. in \cite{synchrosqueeze2}.  Synchrosqueezed wave packets provide more accurate time/space-frequency resolution  representations than classic wavelets by using a reassignment procedure. 
Recently, an arbitrary shape filter bank (ASFB) construction was proposed by Mahde et al. in \cite{asfb1}. Given a boundary line, the ASFB uses optimization techniques to define a low-pass and high-pass filter. However, since the ASFB defines only two filters, it is insufficient for many applications where more modes are present in the image.
One more approach, proposed by Gilles et. al in \cite{1dewt}, \cite{2dewt}, is the empirical wavelet transform (EWT). The empirical wavelet transform aims to build wavelet filter banks whose supports in the frequency domain are detected from the information contained in the spectrum of the signal/image. This process is equivalent to building wavelet filters based on an adaptive partitioning of the Fourier domain. The EWT has multiple bidimensional extensions based on Tensor, Littlewood–Paley, and Curvelet wavelets. The EWT has shown promising results in many applications \cite{ewtunsup},\cite{ewtsup},\cite{ewtapp1},\cite{ewtapp2},\cite{ewtapp3}. However, in 2D, the detected partitions are made of subsets of constrained shapes limiting the adaptability capability of the EWT.  

In this paper, we propose to remove such limitations by defining a new 2D version of the empirical wavelet transform based on fully arbitrary partitions of the frequency domain.
This paper is organized as follows: In section~\ref{sec:2}, we review the empirical wavelet transform and its existing implementations. We illustrate the limitations in the existing 2D approaches. In section~\ref{sec:3}, we investigate the theoretical aspects of building 2D empirical wavelet filters based on arbitrary supports in the Fourier domain. We also propose an algorithm, based on scale-space representations and the watershed transform, to automatically detect a partition corresponding to dominant harmonic modes. Given such a partition, we define the empirical watershed wavelet transform (EWWT). In section~\ref{sec:4}, we illustrate the performance of the EWWT as a multi-resolution analysis method, first on a synthetic image with chosen harmonic modes. Next we illustrate that the EWWT allows improvements on the results obtained in unsupervised texture segmentation and in nonblind deconvolution. Finally, we give our concluding remarks in section~\ref{sec:5}. 

\section{The Empirical Wavelet Transform}
\label{sec:2}

The empirical wavelet transform (EWT) is a method first proposed in 1D by Gilles in \cite{1dewt} and extended to 2D in \cite{2dewt}. The EWT is based on the idea of constructing a family of wavelets based on a set of detected (i.e. data-driven) supports in the Fourier domain corresponding to different harmonic modes. Hereafter, we first recall the 1D empirical wavelet transform, then the 2D extensions, as well as the support detection technique based on scale-space representations.

\subsection{1D Empirical Wavelet Transforms}
\label{sec:2:existing1d}
\begin{figure}[!t]
\centering
\includegraphics[width=.75\textwidth]{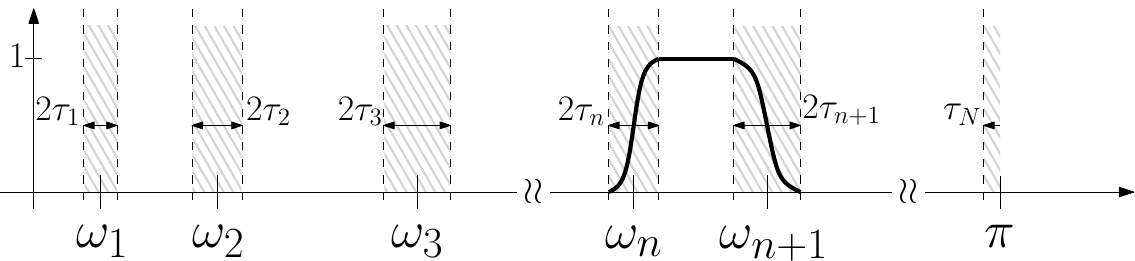}
\caption{Construction of 1D Empirical Wavelet: given a set of boundaries $\{\omega_i$\},  we define $\widehat\psi_n$ as a band-pass filter between $\omega_n$ and $\omega_{n+1}$ with transition regions of width $2\tau_n$ and $2\tau_{n+1}$.}
\label{fig:ewt1d}
\end{figure}

 The one dimensional empirical wavelet transform constructs a set of wavelets based on an arbitrary partitioning of the 1D Fourier domain, see Figure~\ref{fig:ewt1d}. We consider a normalized Fourier axis, which is $2\pi$ periodic, and we denote the Fourier transform of a function $f \in L^1 \cap L^2$ as $\widehat f$. While we consider real signals here, which have symmetric spectra, the generalization to complex signals is fairly straightforward and can be found in \cite{cewt}.

Given a set of detected boundaries $\{\omega_n\}\subset [0,\pi]$ defined on the frequency axis and their corresponding transition widths $\{\tau_n\}$, we define the set of empirical wavelets $\{\widehat\psi_n\}$ as the set of band-pass filters with supports $[-\omega_{n+1},-\omega_n]\cup[\omega_n,\omega_{n+1}]$ and the empirical scaling function $\widehat\phi_1$ as a low-pass filter supported by $[-\omega_1,\omega_1]$. Furthermore, if the transition widths $\tau_n$ are chosen proportionally to $\omega_n$, i.e. $\tau_n = \gamma \omega_n$, where $\gamma < \min_n\Big(\frac{\omega_{n+1}-\omega_n}{\omega_{n+1}+\omega_n}\Big)$, then the obtained set of empirical wavelets $\left\{\phi_1,\{\psi_n\}\right\}$ forms a tight frame of bound 1.

With the set of empirical scaling function and empirical wavelets, $\{\phi_1(t),\{\psi_n(t)\}_{n = 1}^N\}$, we can then define the empirical wavelet transform, denoted $\mathcal{W}^\epsilon_f (n,t)$. The empirical scaling coefficients are given by the inner product of the signal and $\phi_1$, (denoting $h^-(x)=h(-x)$)
\begin{equation}
\label{eq:ewtscale}
\mathcal{W}^\epsilon_f (0,t) = \langle f,\phi_1 \rangle(t) =  (f \ast \phi_1^-)(t) = \Big(\widehat{f} (\omega)\overline{\widehat{\phi_1}(\omega)}\Big)^\lor (t),
\end{equation}
and the detail coefficients $\mathcal{W}^\epsilon_t(n,t)$ are defined similarly,
\begin{equation}
\label{eq:ewtdetail}
\mathcal{W}^\epsilon_f (n,t) = \langle f,\psi_n \rangle(t) =  (f \ast \psi_n^-)(t) = \Big(\widehat{f} (\omega)\overline{\widehat{\psi_n}(\omega)}\Big)^\lor (t),
\end{equation}
where $\ast$ and $( . )^\lor$ denote respectively the convolution operator and the inverse Fourier transform. Due to the tight frame condition, the original signal $f$ can be reconstructed from $\mathcal{W}^\epsilon_f (n,t)$ via
\begin{equation}
\label{eq:ewtreconstruction}
f(t) = \left(\mathcal{W}^\epsilon_f (0,.) \ast \phi_1\right)(t)  + \sum\limits_{n=1}^N\left(\mathcal{W}^\epsilon_f (n,.) \ast\psi_n\right)(t)  = \Big(\widehat{\mathcal{W}^\epsilon_f }(0,\omega)\widehat{\phi}_1(\omega)+ \sum\limits_{n=1}^N \widehat{\mathcal{W}_f^\epsilon}(n,\omega)\widehat{\psi}_n(\omega)\Big)^\lor(t).
\end{equation}

\subsection{2D Empirical Wavelet Transforms}
\label{sec:2:existing2d}
\begin{figure}[!t]
\centering
\subfloat[]{\includegraphics[width=.3\textwidth]{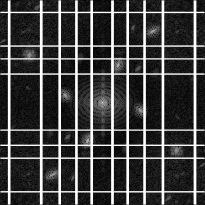}}\hspace{.05in}
\subfloat[]{\includegraphics[width=.3\textwidth]{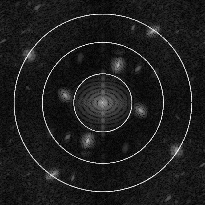}}\hspace{.05in}
\subfloat[]{\includegraphics[width=.3\textwidth]{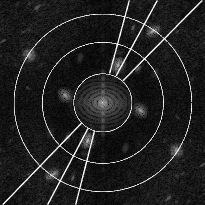}}
\caption{Examples of 2D EWT Implementations: In order from left to right, we see the partitioning of the Fourier spectrum of the \textbf{(a)} Tensor, \textbf{(b)} Littlewood-Paley, and \textbf{(c)} Curvelet (Option 1) empirical wavelet transforms.}
\label{fig:2dexisting}
\end{figure}

In 2D, the set of detected supports is characterized using regions in the Fourier plane, rather than intervals on the real line. Adding this extra dimension brings an extra degree of liberty in the detection of partitions: geometry. In 1D, intervals correspond to basic geometric objects, however in 2D the detected regions can have very different shapes. The 2D extensions proposed in \cite{2dewt} revisit three main families of classic wavelets based on Fourier supports of specific shapes. The first one, the tensor empirical wavelet transform, is a separable approach processing rows and columns independently, based on rectangular supports. The second one aims at building Littlewood-Paley wavelet filters, whose purpose is to separate scales, defined by concentric rings supports. The third one is a curvelet type transform and captures both scale and orientation information by defining supports that are polar wedges shaped supports. The adaptability of these different transforms comes from the fact that the positions and sizes of these different geometric structures are driven by the spectrum of the analyzed image. It is worth noticing that three different strategies were proposed to build the empirical curvelets. The first detects scales and angular sectors independently. The second detects scales first and then detects the angular sectors within each scale ring. The third and final option detects the angular sectors first and then detects scales within each. An example partition of each of the 2D empirical wavelet transforms (only the curvelet first option is depicted here) is illustrated in Figure~\ref{fig:2dexisting}.
If each of these transforms allow for adaptive multi-resolution analysis with specific properties, they are still constrained by the restrictions on the shapes of their detected supports, limiting their degree of adaptability. The purpose of this paper is to remove such limitations by considering completely arbitrary partitions, but before we describe this new approach, the next section will review the scale-space approach to detect such data-driven partitions.

\subsection{Scale-Space Representations}
\label{sec:2:scalespace}
As described in the previous sections, the ability for the EWT to be adaptive lies in the fact that the method is capable of detecting the ``optimal'' supports (in the Fourier domain) for each wavelet filter. This leads to the natural question: how can we detect such a partition of the Fourier domain? Since the goal for each filter is to extract meaningful harmonic modes, let us first define a harmonic mode. A harmonic mode will correspond to an amplitude modulated/frequency modulated (AM/FM) signal, denoted $f_n(t)=A_n(t)\cos(\nu_n(t))$, where $\widehat{A_n}(\omega)$ and $\widehat{\nu_n}(\omega)$ are fast decaying (i.e. of almost bounded support) functions. It is also assumed that two consecutive $\widehat{f_n}(\omega)$ are separated enough to be independently observed. From a geometric point of view, this means that $\left|\widehat{f}(\omega)\right|$ is the superposition of distinct lobes. Detecting the expected partition then corresponds to find the boundaries between successive modes. The original approaches proposed by Gilles et al. in \cite{1dewt} and \cite{2dewt} either consider the midpoint or the lowest minimum between the highest maxima. The major drawback of such approaches is that they require the knowledge of the expected number of modes which is, in the general case, not necessarily known. \\
An automatic detection algorithm (including the number of modes) based on scale-space representations has been proposed in \cite{ewtscalespace} and is recalled hereafter. The scale-space representation \cite{scalespace} of a function $g$ provides the content of $g$ at multiple scales by using a blurring operator to continuously remove components of small scales. In the 1D discrete case, it is given \cite{scalespace} by
\begin{equation}
\label{eq:1dscalespace}
L(t;s) =\sum\limits_{k = -\infty}^\infty T(k;s)g(t-k),
\end{equation}
where $s$ is a scale parameter and $T(k;s) = e^{-\alpha s}I_k(\alpha s)$ is the 1D discrete Gaussian kernel. Note that $I_k$ is the Bessel function of first order and $\alpha \in \mathbb{R}$. A straightforward 2D extension, as will be necessary in Section~\ref{sec:3:modedetection}, is simply obtained by using separable convolutions,
\begin{equation}
\label{eq:scalespace}
L(x,y;s) = \sum\limits_{k = -\infty}^\infty T(k;s) \sum\limits_{l = -\infty}^\infty T(l;s)g(x-l,y-k).
\end{equation}
Since we consider the discrete case, we will sample $s$ with a fixed step-size, denoted $s_0$, which will be discussed in the experimental section. One of main axioms fulfilled by such scale-space representations is the guarantee that the number of extrema can only decrease when $s$ increases. This reduction of extrema permits us to define how meaningful each original minimum is by simply measuring its lifespan (i.e. the value of $s$ before this minimum disappears). Given the lifespans of all initial minima, each of them are then classified as meaningful/persistent or not.

As such, to detect the expected partition of the Fourier domain, we compute the scale-space representation of the signal magnitude spectrum, i.e. we set $g=|\widehat f|$, to extract the meaningful minima, which will correspond to the set of expected boundaries $\{\omega_n\}$. In the 2D extensions described in Section~\ref{sec:2:existing2d}, detecting 2D partitions resumes in performing multiple 1D detections (by averaging the 2D Fourier spectrum or pseudo-polar Fourier spectrum with respect each dimension) because of their specific structures. We invite the reader to check \cite{1dewt}, \cite{2dewt} and especially \cite{ewtscalespace} for more details and numerical implementations.

\section{A New 2D Empirical Wavelet Transform}
\label{sec:3}
In the following sections, we present a general approach to remove the geometric constraints on the detected partitions in the previously described 2D empirical wavelets. First, we give a mathematical construction of wavelet filters based on arbitrarily shaped supports and show that these filters permit us to define a general empirical wavelet transform. We prove the existence of dual filters, allowing the definition of the inverse transform. Second, we present a new algorithm to detect 2D partitions of arbitrary geometry, inspired by the 1D scale-space algorithm discussed in Section~\ref{sec:2:scalespace}.

\subsection{2D empirical wavelets of arbitrary supports}
\label{sec:3:2darbitrary}
Let us assume we are given a discrete image $f\in L^2(\Lambda)$ where $\Lambda=[1,\ldots,\mathcal{N}]\times [1,\ldots,\mathcal{N}]$ is the image domain (pixels coordinates will be denoted $(i,j)\in\Lambda$) and an arbitrary partition, $\{\Omega_n\}_{n=1}^N$, of the Fourier domain (denoted $\Omega=[1,\ldots,\mathcal{N}]\times [1,\ldots,\mathcal{N}]$, frequency coordinates will be denoted $(k,l)\in\Omega$), whose detection will be discussed in Section~\ref{sec:3:ewwt}. 
We assume the partition $\{\Omega_n\}_{n=1}^N$ has the following properties: $\bigcup_{n=1}^N \Omega_n=\Omega$ and $\Omega_n\cap\Omega_m=\emptyset$ if $n\neq m$, whose boundaries delineate the expected supports. In order to define a transition region similar to the one shown in Figure~\ref{fig:ewt1d} for the 1D EWT, we must define the distance from any point in the image to 
these boundaries. As such, for each region $\Omega_n$, we define a distance transform of the region, 
\begin{equation}
\label{eq:distanceTransform}
D_n(k,l) =\begin{cases} 
\frac{2\pi}{\mathcal{N}}\min_{(p,q)\in \partial\Omega_n} \big(d(k,l,p,q)\big) \quad &\text{ if } (k,l) \in \Omega_n \\
-\frac{2\pi}{\mathcal{N}}\min_{(p,q)\in\partial\Omega_n} \big(d(k,l,p,q)\big) &\text{ if } (k,l) \notin \Omega_n
\end{cases},
\end{equation}
where $\partial\Omega_n$ is the boundary of the region $\Omega_n$ and $d(.,.,.,.)$ is the quasi-euclidean distance, defined by
\begin{equation}
\label{eq: qedist}
d(k,l,p,q) = \begin{cases}
(\sqrt2 -1) |q-l| + |p-k| \qquad &\text{ if }\; |p-k| \geq |q - l| \\
(\sqrt2 -1) |p-k| + |q-l|  &\text{ if }\; |p - k|<|q-l|,
\end{cases}
\end{equation}
This is done such that, for each region, we define any point in the image spectrum by its distance to the region's boundary. From there, we can define a 2D empirical wavelet as
\begin{equation}
\label{eq:2dempiricalwavelet}
\widehat\varphi_{n,k,l} = \begin{cases} 
1  &\text{ if }  D_n(k,l) > \tau \\ 
\cos\left( \frac{\pi}{2} \beta\left( \frac{\tau-  D_n(k,l)}{2\tau} \right)\right) \quad &\text{ if } D_n(k,l) \leq |\tau| \\ 
0 &\text{ if }   D_n(k,l) < -\tau,
\end{cases}
\end{equation}
where $\tau$ is the desired transition area width and $\beta(x)$ is an arbitrary $\mathcal{C}^k$ function such that
\begin{equation}
\label{eq:betax}
\beta(x) = \begin{cases}
0 \hspace{.05in} \text{ if } x \leq 0 \\
1 \hspace{.05in} \text{ if } x\geq 1 \end{cases} \quad\text{and}\quad \beta(x) + \beta(1-x) = 1 \hspace{.05in}\quad \forall x \in [0,1] .
\end{equation}
A usual choice for $\beta$ is
\begin{equation}
\label{eq:beta2}
\beta(x) = x^4(35-84x+70x^2-20x^3).
\end{equation}
 Unlike in the 1D case or the existing 2D cases, there are no theoretical results for the appropriate choice of $\tau$, and therefore it must be chosen experimentally. Nevertheless, we can show that the set $\{\varphi_n\}$ forms a frame.

\begin{Proposition}
Denoting $\varphi_{n,i,j}=\varphi_n(\cdot-i,\cdot-j)$, the set $\{\varphi_{n,i,j}\}$ forms a frame.
\end{Proposition}
\begin{proof}
Assuming, for now, that (to be shown later) $\exists A,B\in\mathbb{R},0< A \leq B < \infty$, such that $\forall (k,l)\in\Omega$
$$A\leq \sum_{n=1}^N |\widehat\varphi_{n,k,l}|^2 \leq B,$$
then, using Parseval's theorem, we get 
\begin{gather}
 A|\widehat{f}_{k,l}|^2\leq \sum_{n=1}^N \left|\widehat{f}_{k,l}\overline{\widehat{\varphi}_{n,k,l}}\right|^2\leq B|\widehat{f}_{k,l}|^2\\
\Rightarrow A\sum_{(k,l)\in\Omega}|\widehat{f}_{k,l}|^2\leq \sum_{n=1}^N \sum_{(k,l)\in\Omega}\left|\widehat{f}_{k,l}\overline{\widehat{\varphi}_{n,k,l}}\right|^2\leq B\sum_{(k,l)\in\Omega}|\widehat{f}_{k,l}|^2\\
\Rightarrow A\|f\|^2\leq \sum_{n=1}^N \sum_{(i,j)\in\Lambda}\left|\left\langle f, \varphi_{n,i,j}\right\rangle\right|^2\leq B \|f\|^2
\end{gather}
which implies that $\{\varphi_{n,i,j}\}$ forms a frame of bounds $A$ and $B$.\\
It remains to show that such bounds $A$ and $B$ exist such that $\forall (k,l)\in\Omega,A\leq \sum_{n=1}^N |\widehat \varphi_{n,k,l}|^2\leq B$.
We can observe that three situations, as shown in Figure~\ref{fig:proofcases}, can occur:
\begin{itemize}
\item[-] $(k,l)$ lies outside of a transition area
\item[-] $(k,l)$ lies in a transition area between only two regions
\item[-] $(k,l)$ lies in a transition area between three or more regions
\end{itemize}
\begin{figure}[!t]
\centering
\includegraphics[width=1.75in]{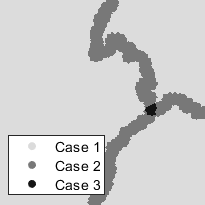}
\caption{Three Cases Described In Proof: Depending on whether $(k,l)$ falls within a transition region, and how many regions are in that transition, we see different bounds for our frame.}
\label{fig:proofcases}
\end{figure}
It is easy to see that if $(k,l)$ is not within a transition area, $\varphi_{n,k,l}$ is nonzero for only one index $n$, and by construction, we have $\sum_{n=1}^N |\widehat\varphi_{n,k,l}|^2 = 1$. It remains to show the bounds if $(k,l)$ lies within a transition area.

The lower bound is decided when the transition area is defined between two regions. In this situation, we have 
$$|\widehat\varphi_{1,k,l}|^2 + |\widehat\varphi_{2,k,l}|^2 = \Big[\cos\Big( \frac{\pi}{2} \beta\Big( \frac{\tau-  D_1(k,l)}{2\tau}\Big)\Big)\Big]^2 + \Big[\cos\Big( \frac{\pi}{2} \beta\Big( \frac{\tau-  D_2(k,l)}{2\tau}\Big)\Big)\Big]^2.$$ 
But, from the definition of the distance transform, we have $D_2(k,l) = -D_1(k,l)$, hence
$$|\widehat\varphi_{1,k,l}|^2 + |\widehat\varphi_{2,k,l}|^2 = \Big[\cos\Big( \frac{\pi}{2} \beta\Big( \frac{\tau-  D_1(k,l)}{2\tau}\Big)\Big)\Big]^2 + \Big[\cos\Big( \frac{\pi}{2} \beta\Big( \frac{\tau+  D_1(k,l)}{2\tau}\Big)\Big)\Big]^2.$$
Since $\frac{\tau+  D_1(k,l)}{2\tau} + \frac{\tau-  D_1(k,l)}{2\tau} = 1$ and $\beta(x) + \beta(1-x) = 1$, 
$$|\widehat\varphi_{1,k,l}|^2 + |\widehat\varphi_{2,k,l}|^2 =
\Big[\cos\Big( \frac{\pi}{2} \beta\Big( \frac{\tau-  D_1(k,l)}{2\tau}\Big)\Big)\Big]^2 + \Big[\cos\Big( \frac{\pi}{2} \Big(1 - \beta\Big( \frac{\tau-  D_1(k,l)}{2\tau}\Big)\Big)\Big)\Big]^2,$$
and since $\cos(x) = \sin(x+\frac{\pi}{2})$, this is equivalent to 
$$|\widehat\varphi_{1,k,l}|^2 + |\widehat\varphi_{2,k,l}|^2 = \Big[\cos\Big( \frac{\pi}{2} \beta\Big( \frac{\tau-  D_1(k,l)}{2\tau}\Big)\Big)\Big]^2 +\Big[\sin\Big( \frac{\pi}{2}\beta\Big( \frac{\tau-  D_1(k,l)}{2\tau}\Big)\Big)\Big]^2 = 1.$$
Therefore, our lower bound is $A = 1$. 

In the situation where the transition area is between $r$ regions, we have  $|\widehat\varphi_{1,k,l}|^2 + |\widehat\varphi_{2,k,l}|^2 + |\widehat\varphi_{3,k,l}|^2 + \dots + |\widehat\varphi_{r,k,l}|^2 = 1 + |\widehat\varphi_{3,k,l}|^2+\dots + |\widehat\varphi_{r,k,l}|^2$. Since $\forall n = 1, 2, \dots, N,\forall (k,l)\in\Omega,|\widehat\varphi_{n,k,l}|\leq 1$, 
we therefore get an upper bound $B = r-1$ which completes the proof.
\end{proof}

With our set of empirical wavelets in hand, the corresponding transform is then defined by 
\begin{equation}
\mathcal{W}^\epsilon_f (n,i,j) = \langle f,\varphi_{n,i,j} \rangle = (f \ast \varphi_{n,\cdot-i,\cdot-j}^-) = \Big(\widehat{f}_{k,l}\overline{\widehat{\varphi}_{n,k,l}}\Big)^\lor(i,j).
\end{equation}
Since $\forall(k,l)\in\Omega, \sum_{n=1}^N |\widehat\varphi_{n,k,l}|^2>0$, we can define the corresponding dual frame $\{\tilde{\varphi}_{n,i,j}\}$ via
$$\widehat{\tilde\varphi}_{n,k,l} = \frac{\widehat\varphi_{n,k,l}}{\sum_n |\widehat\varphi_{n,k,l}|^2}.$$
The existence of such dual frame allows us to reconstruct the original signal $f$ from its empirical wavelet transform:
\begin{equation}
f_{i,j} =\sum_{n=1}^N \left(\mathcal{W}_f^\epsilon(n,\cdot,\cdot)\ast \tilde\varphi_{n,\cdot-i,\cdot-j}\right) = \Big(\sum_{n=1}^N  \widehat{\mathcal{W}}^\epsilon_f(n,k,l) \widehat{\tilde\varphi}_{n,k,l}\Big)^\lor(i,j).
\end {equation}

\subsection{On the detection of partitions of arbitrary geometry}
\label{sec:3:ewwt}
In this section, we address the question of how to detect partitions of the Fourier domain corresponding to meaningful harmonic modes of arbitrary geometry. The proposed method uses a 2D scale-space representation to find the ``center'' of harmonic modes; followed by a watershed transform to find the arbitrary boundaries delimiting the different regions $\Omega_n$. These two steps are described in detail in the next two sections. From now on, we assume our image $f$ is discrete of size $\mathcal{N} \times \mathcal{N}$.

\subsubsection{Scale-space localization of harmonic modes}
\label{sec:3:modedetection}
\begin{figure}[!t]
\centering
\subfloat[]{\includegraphics[width=.49\textwidth]{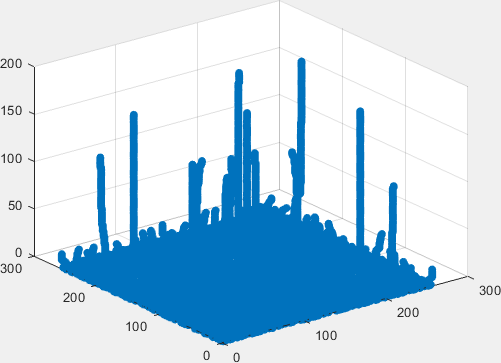}}\hspace{.05in}
\subfloat[]{\includegraphics[width=.49\textwidth]{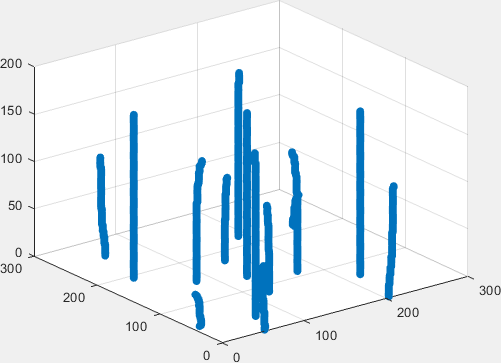}}
\caption{Tracking of Maxima through Scale-Space (the vertical axis corresponds to the scale $s$): \textbf{(a)} shows all maxima through scale-space. \textbf{(b)} shows the remaining persistent maxima, i.e of lifespan larger than the threshold defined by Otsu's method.}
\label{fig:maximatracking}
\end{figure}

In Section~\ref{sec:2:scalespace}, we showed that the expected 1D boundaries corresponded to meaningful local minima which can be found using scale-space representations. Unfortunately, 2D boundaries separating the different regions correspond to arbitrary curves and it becomes very hard to characterize what a meaningful curve is. Since our goal is to separate harmonic modes, assuming that these modes are well enough separated lobes in the spectrum, we propose first to detect the potential candidates by selecting meaningful local maxima. To do so, we follow a similar approach as in Section~\ref{sec:2:scalespace} but looking for persistent maxima instead of minima. To resume the process, we build a 2D scale-space representation of the image magnitude spectrum, i.e. $g=|\widehat f|$, detect all local maxima at each scale $s$, and measure their lifespans (i.e the largest $s$ before they disappear). From there, we create a histogram of the persistence of the maxima and use Otsu's method to define a threshold to classify each maxima as persistent or not. The locations of all persistent maxima, denoted $\{\mu_n\}$, are then extracted; these locations represent ``centers'' of the expected harmonic modes. Because we are working in a discrete space, the scale step-size $s_0$ is a parameter and we set the maximum scale proportional to $s_0$ and the image size, $s_{max} = 4\frac{s_0\mathcal{N}}{K}$, where $K$ is the size of the kernel used to create the scale-space representation.
Figure~\ref{fig:maximatracking}.a illustrates the tracking of maxima through the scale-space, while Figure~\ref{fig:maximatracking}.b shows the remaining persistent maxima after thresholding. 

\subsubsection{Watershed Transforms}
\label{sec:3:watershed}
Given a set of mode centers, $\{\mu_n\}$, we wish to find a set of boundaries which will define the mode supports. Taking again the point of view that modes correspond to lobes, then a natural way to find such boundaries is to select the bottom of the valleys between the modes. From a mathematical perspective, such a process corresponds to finding the path of lowest separation (this idea is a direct extension of the principle used in \cite{2dewt}, for finding 1D boundaries, where the lowest minimum between meaningful peaks was chosen). 

The watershed transform, first proposed by Beucher et al. in \cite{watershed1}, is an image segmentation technique that defines a contour based on the path of highest separation between minima. Based on the geographic definition of watersheds and catchment basins, the watershed transform treats an image as a topographic landscape with pixel intensity representing the height at that pixel. Then, the transforms separates the image into its catchment basins, with a watershed contour separating them. This process can be described as follows.

Given an image $g$ and a set of its minima, $\{x_i\}$, we uniformly ``flood'' the topographic landscape  produced from $g$, with the water collecting at the minima. When one body of water would flow into another, we construct a barrier. Once complete, we define the contour $\Gamma$ along the barrier.

Later, Meyer et al. in \cite{watershed2} proposed a method in which one morphologically reconstructs the image in such a way to impose minima at select markers $\{M_n\}$. This is done by first forcing each marker to be a minimum, by setting $g(M_n)= -\infty$ for all $n$, and then by filling the catchment basins of each minimum $x_i\notin \{M_n\}$ in such a way that the region contains no extrema. This approach reconstructs the image so that the only minima are at the markers $\{M_n\}$ and thus reformulates the watershed transform such that the generated contour $\Gamma$ lays along the path of highest separation between selected markers, rather than all minima. 

To find the set of boundaries that will define the expected supports, we first invert the magnitude spectrum, $f_- = -|\widehat f|$, so that the watershed transform will define the path of lowest separation, rather than highest separation. Then, we impose minima to be at the location of persistent maxima, $\{\mu_n\}$. Finally, we apply the watershed transform, which defines a boundary $\Gamma$ on the magnitude spectrum $|\widehat f|$ that is along the paths of lowest separation between the locations of persistent maxima $\{\mu_n\}$. The boundary $\Gamma$ defines a partition with regions $\{\Omega_n\}$.

Since the watershed transform defines some pixels as part of the boundary, we must assign these pixels to a region. This assignment is not critical, since they belong to the path of lowest separation, but it should be symmetric if $f$ is a real valued image. Furthermore, if $f$ is real valued, we must pair regions symmetrically with respect to the origin. This can be achieved using Algorithm~\ref{alg:pairing}.

\begin{algorithm}[H]
\begin{algorithmic}
\STATE \textbf{Input: } A set of unpaired maxima locations $\{\mu_n\}$ and corresponding partition $\{\Omega_n\}$, where $n = 1, 2, \dots,N$ \\[.1in]
\STATE $m_P \leftarrow \begin{cases} \frac{N+1}{2} \text{ if } (0,0) \in \{\mu_n\} \\[.1in]
	\frac{N}{2} \text{ otherwise}\end{cases}$\\[.1in]
\FOR{$n = 1,2,\dots N$}
	\STATE $\Theta_n = \{m | \pm\mu_n \in  \Omega_m\}$\\[.1in]
	\STATE $\widetilde \Omega_n = \bigcup\limits_{m\in\Theta_n} \Omega_m$\\[.1in]
	\STATE $\{\mu_n\} \leftarrow \{\mu_n\} - \{\mu_{\Theta_n}\}$\\[.1in]
\ENDFOR \\[.1in]
\STATE \textbf{Output:} Set of paired regions $\{\Omega_n\} \leftarrow \{\widetilde \Omega_n\}$, where $n = 1,2,\dots,N$
\end{algorithmic}
\caption{Boundary pixel assignment and regions symmetrization}
\label{alg:pairing}
\end{algorithm}

\subsection{The Empirical Watershed Wavelet Transform}
\label{sec:3:ewwtfinal}
\begin{figure}[!t]
\centering
\subfloat[]{\includegraphics[width = .11\textwidth]{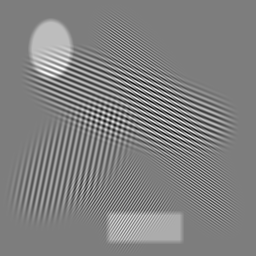}\label{sbsa}}\hspace{.05in}
\subfloat[]{\includegraphics[width = .11\textwidth]{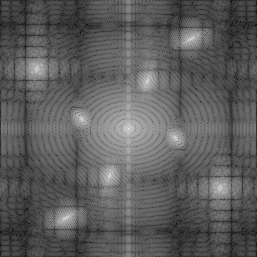}\label{sbsb}}\hspace{.05in}
\subfloat[]{\includegraphics[width = .11\textwidth]{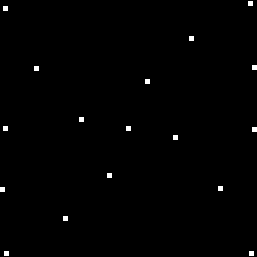}\label{sbsc}}\hspace{.05in}
\subfloat[]{\includegraphics[width = .11\textwidth]{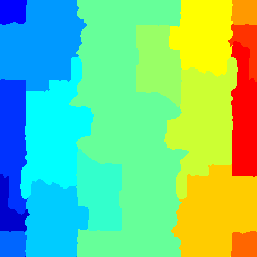}\label{sbsd}}\hspace{.05in}
\subfloat[]{\includegraphics[width = .11\textwidth]{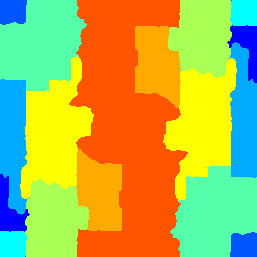}\label{sbse}}\hspace{.05in}
\subfloat[]{\includegraphics[width = .11\textwidth]{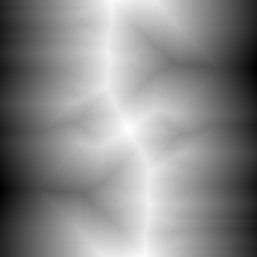}\label{sbsf}}\hspace{.05in}
\subfloat[]{\includegraphics[width = .11\textwidth]{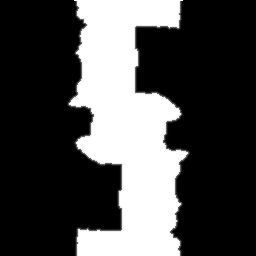}\label{sbsg}}\hspace{.05in}
\subfloat[]{\includegraphics[width = .11\textwidth]{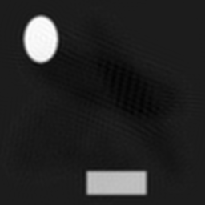}\label{sbsh}}
\caption{Step by step construction of empirical watershed wavelets: \textbf{(a)} Original image, \textbf{(b)} Magnitude spectrum of image, \textbf{(c)} Detected persistent maxima using method from Section~\ref{sec:3:modedetection}, \textbf{(d)} detected partitioning of the spectrum using method from Section~\ref{sec:3:watershed}, \textbf{(e)} Paired regions if image is real, \textbf{(f)} Distance transform of central region, \textbf{(g)} Empirical watershed wavelet of central region, \textbf{(h)} empirical watershed wavelet coefficient for central region. \label{fig:fullprocess}}
\end{figure}

We are now equipped with all tools to define a new 2D arbitrary empirical wavelet transform called the empirical watershed wavelet transform (EWWT). The process is described by Algorithm~\ref{alg:ewwt}, and the output of the intermediate steps are illustrated in Figure~\ref{fig:fullprocess}.

\begin{algorithm}[H]
\begin{algorithmic}
\STATE \textbf{Input: } Image $f$ (Figure~\ref{sbsa}), transition thickness $\tau$ and scale-space step-size $s_0$\\[.1in]
\STATE Use scale-space approach with step-size $s_0$ from Section \ref{sec:3:modedetection} to find the location of persistent maxima $\{\mu_n\}$ (Figure~\ref{sbsc})   \\[.1in]
\STATE Use watershed transform from Section \ref{sec:3:watershed} to detect the partition $\{\Omega_n\}$  (Figure~\ref{sbsd})\\[.1in]
\IF{$f$ is real valued image} 
    \STATE{Pair $\{\Omega_n\}$ using the pairing algorithm in Section \ref{sec:3:watershed}}  (Figure~\ref{sbse})\\[.1in] 
    \ENDIF \\[.1in]
\STATE Build the wavelet filters $\{\widehat{\varphi_n}\}$ of supports $\{\Omega_n\}$ using the construction described in \ref{sec:3:2darbitrary} (Figure~\ref{sbsg})\\[.1in]
\STATE Compute the EWWT given by $\mathcal{W}^\epsilon_f (n,i,j) = \Big(\widehat{f}_{k,l}\overline{\widehat{\varphi}_{n,k,l}}\Big)^\lor (i,j)$ (Figure~\ref{sbsh})\\[.1in]
\STATE \textbf{Output: } Set of empirical watershed wavelet coeficients $\{\mathcal{W}^\epsilon_f(n,i,j)\}$ and its corresponding partition $\{\Omega_n\}$ and filter bank $\{\widehat{\varphi_n}\}$.
\end{algorithmic}
\caption{The Empirical Watershed Wavelet Transform.}
\label{alg:ewwt}
\end{algorithm}

\section{Experiments}
\label{sec:4}
In this section, we first qualitatively explore the performance of the EWWT by decomposing a toy image made of known harmonic modes. Next, we propose to quantitatively test the use of the EWWT on two practical applications: unsupervised texture segmentation and non-blind deconvolution. As stated in Section~\ref{sec:3:2darbitrary} and \ref{sec:3:modedetection}, the transition width $\tau$ and scale-space step-size $s_0$ must be provided. Unless otherwise specified, we choose $\tau = 0.1$ and $s_0 = 1$ for all experiments. 

\subsection{Harmonic mode decomposition}

\begin{figure}[!t]
\centering-
\subfloat[]{\includegraphics[width=.25\textwidth]{Figures/origImg.png}\label{fig:origimg}}\hspace{.05in}
\subfloat[]{\includegraphics[width=.25\textwidth]{Figures/Spectrum.png}\label{fig:origspectrum}}\hspace{.05in}
\subfloat[]{\includegraphics[width=.25\textwidth]{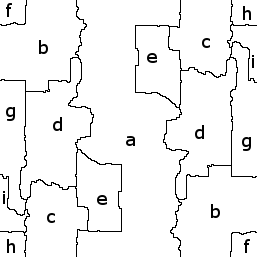}\label{fig:partition}}
\caption{\textbf{(a)} An example image, \textbf{(b)} its magnitude Fourier transform, \textbf{(c)} the detected partition using the EWWT.}
\end{figure}
\begin{figure}[!t]
\centering
\subfloat[$L^2 = 129.8058$]{\includegraphics[width=.25\textwidth]{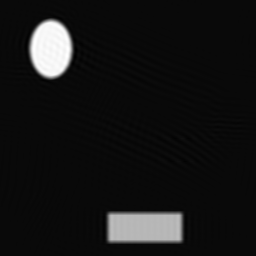}\label{ewwta}}\hspace{.05in}
\subfloat[$L^2 = 7.4679$]{\includegraphics[width=.25\textwidth]{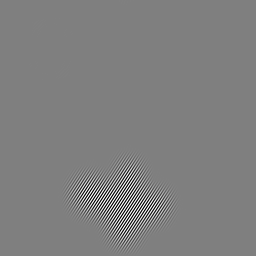}\label{ewwtb}}\hspace{.05in}
\subfloat[$L^2 = 14.135$]{\includegraphics[width=.25\textwidth]{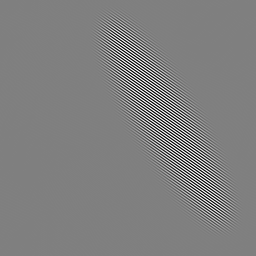}\label{ewwtc}}\\
\subfloat[$L^2 = 9.58$]{\includegraphics[width=.25\textwidth]{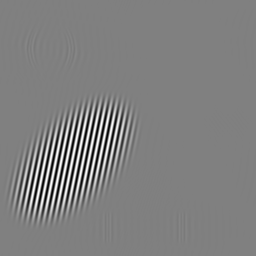}\label{ewwtd}}\hspace{.05in}
\subfloat[$L^2 = 18.380$]{\includegraphics[width=.25\textwidth]{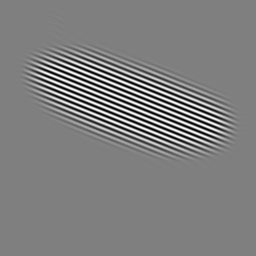}\label{ewwte}}\hspace{.05in}
\subfloat[$L^2 = 0.0019$]{\includegraphics[width=.25\textwidth]{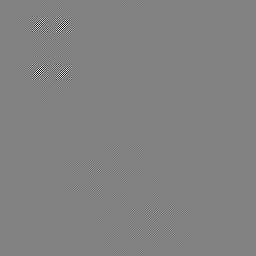}\label{ewwtf}}\\
\subfloat[$L^2 = 0.0511$]{\includegraphics[width= .25\textwidth]{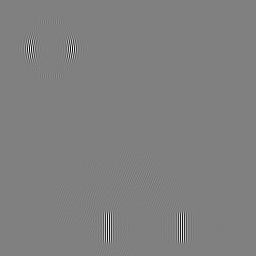}\label{ewwtg}}\hspace{.05in}
\subfloat[$L^2 = 0.0018$]{\includegraphics[width=.25\textwidth]{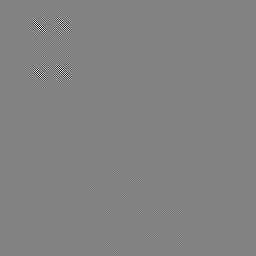}\label{ewwth}}\hspace{.05in}
\subfloat[$L^2 = 0.0074$]{\includegraphics[width=.25\textwidth]{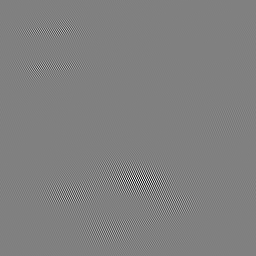}\label{ewwti}}\hspace{.05in}
\caption{Empirical watershed wavelet coefficients for the image shown in Figure~\ref{fig:origimg} with the partitioning shown in Figure~\ref{fig:partition}. These images have been normalized for visibility, but their $L^2$ energy is given to appreciate how much the corresponding information is significant.}
\label{fig:EWWTex}
\end{figure}

First, we show how the EWWT decomposes a toy image, which is made of four well separated harmonic modes and two slightly blurred piecewise constant components. This toy image is illustrated in Figure~\ref{fig:origimg} along with its magnitude spectrum in Figure~\ref{fig:origspectrum}. Figure~\ref{fig:partition} provides the detected partition and Figure~\ref{fig:EWWTex} shows the empirical watershed wavelet coefficients obtained by performing the transform. As expected, the EWWT is capable of isolating the five harmonic modes present in the original image, which are depicted in Figure~\ref{fig:EWWTex}\subref{ewwta}--\subref{ewwte}. As for the other coefficients shown in Figure~\ref{fig:EWWTex}\subref{ewwtf}--\subref{ewwti}, we see very little information (as per their respective $L^2$ norms), which mostly consists of high frequency information near the piecewise regions. Looking at the Fourier partition, we see that these extra coefficients correspond to regions near the edge of the spectrum where the Fourier coefficients are very small.

\subsection{Performance on Unsupervised Texture Segmentation}
In this section, we explore the use of the EWWT to extract texture features for the task of unsupervised texture segmentation. Wavelets have been widely used in the literature to characterize textures since they can extract multi-scale as well as directional information. Recently, in \cite{ewtunsup}, Huang et al. showed the effectiveness of the empirical wavelet transform on this problem, particularly that the Curvelet-1 EWT option outperforms all traditional wavelets. As such, we compare the EWWT against the Curvelet-1 option on the Outex dataset \cite{outex}. 

To perform the unsupervised texture segmentation of an image, we go through the following steps, which are outlined in detail in \cite{ewtunsup}. First, we perform a cartoon-texture decomposition of the image and keep only the textural part. This removes some variability in lighting and other non-textural information. Next, we perform the empirical wavelet feature extraction, which in our case, is either the EWWT or the EWT-Curvelet 1 option. Afterwards, we perform an energy aggregation method by calculating the $L^2$ energy on a window centered at each pixel. The ideal window size is mostly dependent on the dataset, and \cite{ewtunsup} shows that a $19\times19$ window is best for the Outex dataset. The final segmentation is obtained by performing a pixel-wise $k$-means (the expected number of textures $k$ is provided to the algorithm) clustering with \emph{cityblock} distance. 

To quantitatively evaluate the quality of the segmentation, we perform the segmentation on 100 images and compute the score, $R$, which averages six different region-based segmentation benchmarks (the closer to 100\% the better the segmentation). These metrics are the normalized variation of information \cite{nvoi}, the swapped directional hamming distance \cite{sdhd}, the van Dongen distance \cite{vd}, the swapped segmentation covering \cite{ssc1} \cite{ssc2}, the bipartite graph matching \cite{bgm1} \cite{bgm2}, and the bidirectional consistency error \cite{bce}. Each of these metrics are implemented by the \emph{eval\_segm} function available in the Supervised Evaluation of Image Segmentation Methods (SEISM) toolbox\footnote{https://github.com/jponttuset/seism}. 

\subsubsection{Influence of the EWWT parameters}
\label{ch4:ewwtparams}
First, we investigate the choice of the scale-space step-size $s_0$. When using smaller values for $s_0$, we automatically reduce the maximum scale, $s_{max}$. In these cases, Otsu's method detects an earlier threshold, providing a larger number of modes. Table~\ref{tab:earlyscales} shows the results for different values of $s_0$. We observe that a lower $s_0$ outperforms a larger $s_0$, meaning that detecting a large number of modes is advantageous. This suggests that the empirical watershed wavelet transform does a successful job of well separating relevant modes that characterize the textural information. Ultimately, a value of $s_0 = 0.1$ is the most successful, and we will use this choice going forward. 

\begin{table}[H]
 
\caption{Influence of truncating scale-space step-size $s_0$ for texture segmentation (best score in bold)}.\label{tab:earlyscales}
\centering
\begin{tabular}{cccccccccccccc}
\toprule
$\boldsymbol{ s_0}$ & \textbf{0.05} & \textbf{0.1} & \textbf{0.2} &\textbf{0.3} & \textbf{0.4} & \textbf{0.5}  &\textbf{0.6} &\textbf{0.7} & \textbf{0.8} & \textbf{0.9} & \textbf{1.0}\\
\midrule
$R$ & 87.22 &\textbf{87.95} & 87.30 & 86.09 & 86.09 & 85.32  & 85.17 & 83.54 & 80.24 & 76.46 & 73.44 & \\
\bottomrule
\end{tabular}
\end{table}

The next parameter influence we consider is the choice of the transition width $\tau$. Table~\ref{tab:tau} shows that a transition width of $\tau = 0.1$ performs best, which is likely a balance between proper separation of modes and the forming of Gibbs phenomena. 
\begin{table}[H]
\caption{Influence of transition width $\tau$ for texture segmentation (best score in bold)}.\label{tab:tau}
\centering
\begin{tabular}{cccccccccccc}
\toprule
$\boldsymbol \tau$ & \textbf{0.05} & \textbf{.1} & \textbf{.2} & \textbf{.3} & \textbf{.5}\\
\midrule
$R$  & 87.51 & \textbf{87.95} & 87.78 & 86.42 & 87.21 \\
\bottomrule
\end{tabular}
\end{table}

We also look into the effect of excluding maxima detected at the Fourier domain edge in the mode detection step. If a mode is detected at the edge of the spectrum, it represents a mode with very high frequencies. Unless equally high frequency textures are present in the dataset, it is likely these represent noise rather than meaningful textural information. Therefore, we propose to remove the locations $\mu_n$ corresponding to maxima within a certain band around the edge of the Fourier domain. Setting $s_0 = 0.1$ and $\tau = 0.1$, Table~\ref{tab:edges} investigates the influence of the width of that exclusion band. We can observe that removing the maxima within a two pixels distance from the edge of the Fourier domain gives the best accuracy.
	
\begin{table}[H]
\caption{Influence of edge exclusion width for texture segmentation (best score in bold)}.\label{tab:edges}
\centering
\begin{tabular}{cccccccc}
\toprule
\textbf{Pixels}& \textbf{0} & \textbf{1} & \textbf{2} &\textbf{3} &\textbf{ 5}\\
\midrule
$R$ & 87.95 & 88.00 & \textbf{88.37} & 87.16 & 87.31 \\
\bottomrule
\end{tabular}
\end{table}

\subsubsection{Final comparison}
\begin{figure}[!t]
\centering
\subfloat[Original]{\includegraphics[width=.2\textwidth]{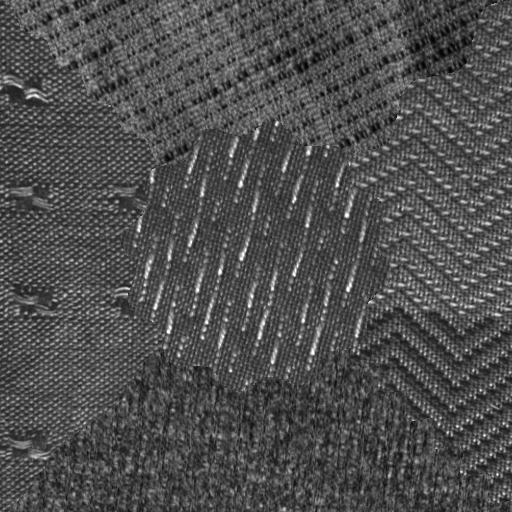}\label{a}} 
\subfloat[Ground Truth]{\includegraphics[width=.2\textwidth]{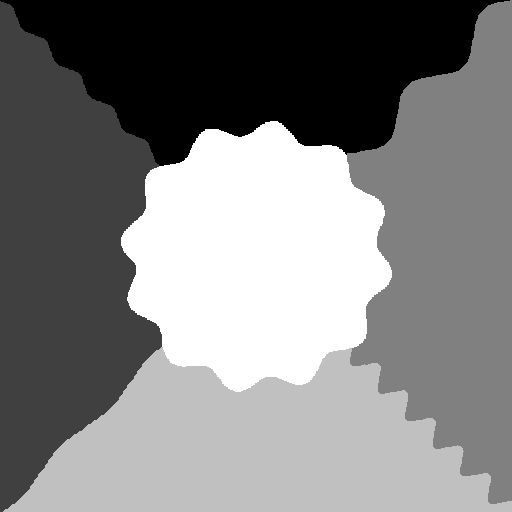}\label{b}}
\subfloat[EWT-C1]{\includegraphics[width=.2\textwidth]{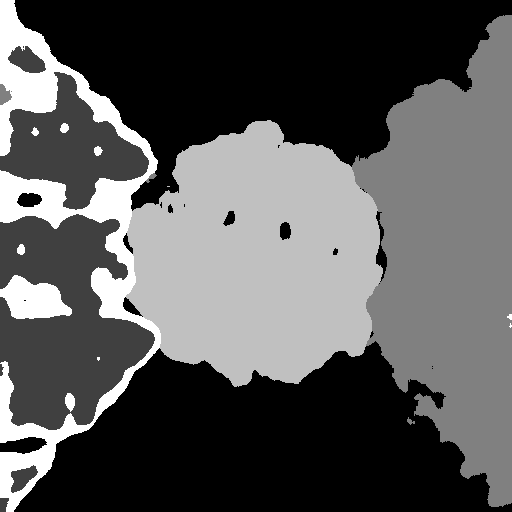}\label{c}} 
\subfloat[EWWT]{\includegraphics[width=.2\textwidth]{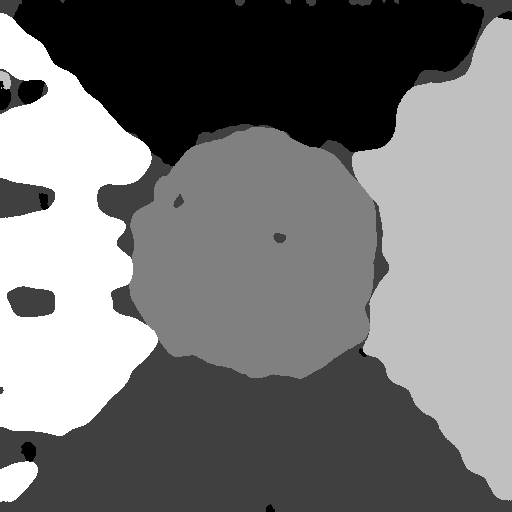}\label{d}}\\
\subfloat[Original]{\includegraphics[width=.2\textwidth]{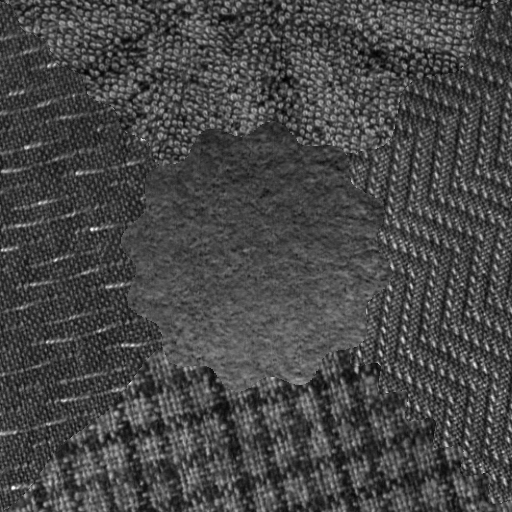}\label{e}}
\subfloat[Ground Truth]{\includegraphics[width=.2\textwidth]{Figures/segex/GT.png}\label{f}}
\subfloat[EWT-C1]{\includegraphics[width=.2\textwidth]{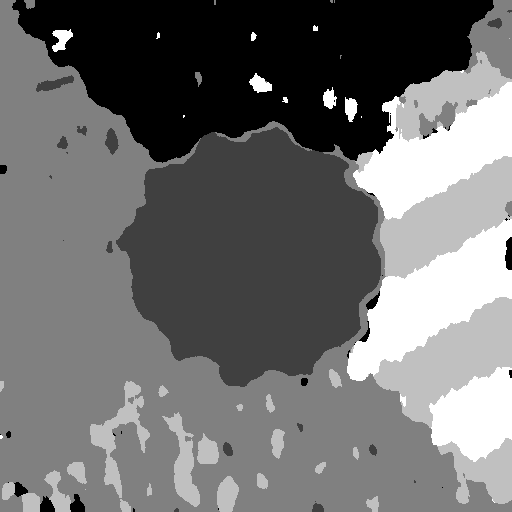}\label{g}}
\subfloat[EWWT]{\includegraphics[width=.2\textwidth]{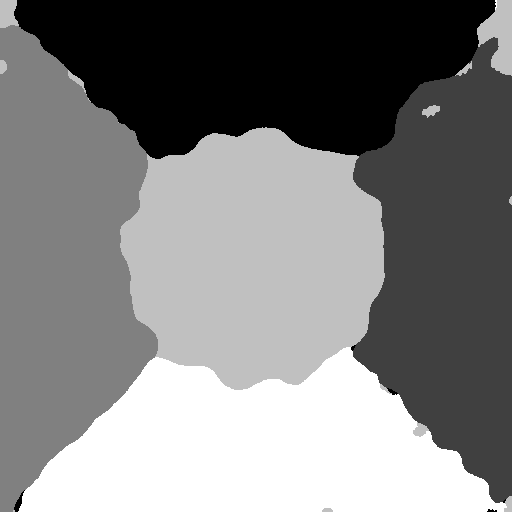}\label{h}}\\
\subfloat[Original]{\includegraphics[width=.2\textwidth]{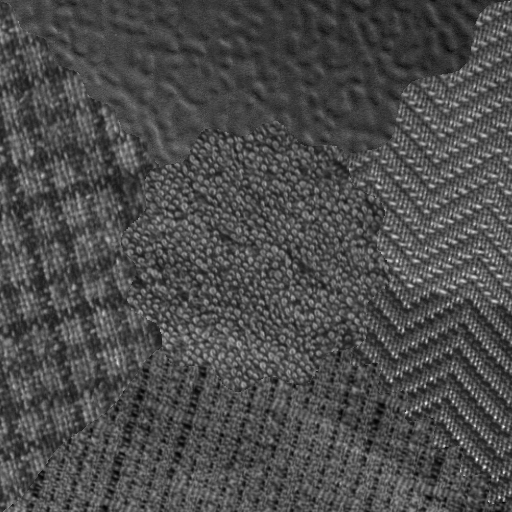}\label{i}}
\subfloat[Ground Truth]{\includegraphics[width=.2\textwidth]{Figures/segex/GT.png}\label{j}}
\subfloat[EWT-C1]{\includegraphics[width=.2\textwidth]{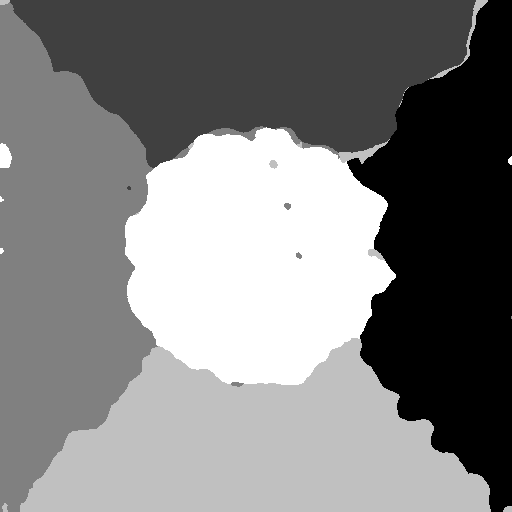}\label{k}}
\subfloat[EWWT]{\includegraphics[width=.2\textwidth]{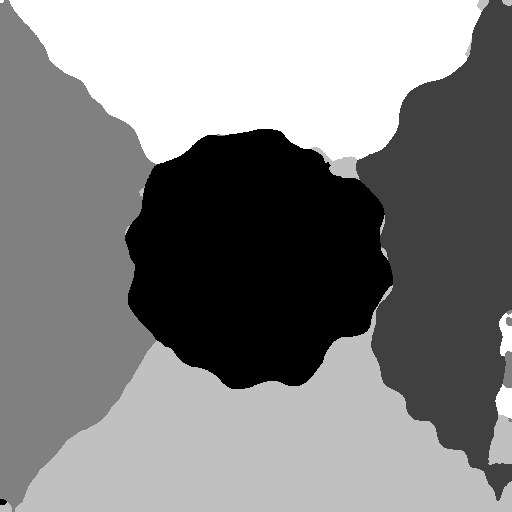}\label{l}}\\
\caption{Results of Segmentation performed on Outex using empirical curvelets (EWT-C1) and empirical watershed wavelets (EWWT).}
\label{fig:segmentEx}
\end{figure}

We now compare the EWWT against the Curvelet-1 EWT (extensive comparisons with other types of wavelets are provided in \cite{ewtunsup}) on the Outex dataset. We use the best parameters for the EWWT for these tests, i.e. we set $s_0 = 0.1$, $\tau = 0.1$, and we exclude maxima detected within two pixels of 
the domain edge. The results are shown in Table~\ref{tab:final}, and we see a significant improvement (almost 2\% with similar standard deviation) when using the empirical watershed wavelet transform. Examples of obtained segmentations can be seen in Figure~\ref{fig:segmentEx}.
\begin{table}[H]
\caption{Final comparative results for texture segmentation (best score in bold)}.\label{tab:final}
\centering
\begin{tabular}{cccc}
\toprule
\textbf{Method} & \textbf{Curvelet-1} & \textbf{EWWT}\\
\midrule
$R$ & 86.41 (7.09) & \textbf{88.37 (7.50)} \\
\bottomrule
\end{tabular}
\end{table}
Likely, this improvement is due to extra adaptability of the EWWT in both mode detection and boundary detection. 

\subsection{Performance of EWWT for Non-Blind Deconvolution}
\label{ch4sec3}
Another standard application of traditional wavelets is in deconvolution as a sparsity-promoting operator\cite{cai1}. Often, an ideal image $u$ is affected by a kernel $A$ and noise $\mu$, which is observed as the image $f$ accordingly to $f = A \ast u + \mu$, where $\ast$ denotes convolution. As such, deconvolution is an inverse problem for image restoration where the goal is, from the observed image $f$, to recover the unknown pristine image $u$. Here we consider the case of non-blind deconvolution where the convolution kernel $A$ is assumed to be known. We begin by discussing how deconvolution inverse problems are solved using traditional wavelets.

\subsubsection{Deconvolution with Traditional Wavelets}
When attempting to recover a pristine image $u$ via deconvolution, we are essentially solving $\min\limits_{u} ||Au-f||_2$. However, this equation is ill-posed, as it lacks a regularization term. Cai et al \cite{cai1} claim that real images usually have sparse representations under wavelet transforms. Therefore, by including a regularization term of $||D u||_1$, where $D$ is a traditional wavelet transform, we ensure that the resulting $u$ will have the characteristics of a real image. The goal is to ultimately turn the image deconvolution problem into a saddle point problem so that we may employ the Primal-Dual method \cite{primaldual}. As such, we also employ a quadratic penalty method, leaving us with 
\begin{equation}
\label{eq:deconvolution}
\min\limits_u ||Du||_1 + \frac{\lambda}{2}||A\ast u-f||_2^2,
\end{equation}
where $\lambda$ is a regularization parameter.

For the reader's convenience, we recall that the Primal-Dual algorithm \cite{primaldual} is a technique to solve saddle point problems of the form 
\begin{equation}\min_u F(Ku) + G(u),\end{equation}
where $K:X\to Y$ is a linear operator. The minimizer $u$ is then obtained by Algorithm~\ref{alg:pd} (where $\operatorname{prox}_{\sigma F}(x) = \underset{y}{\operatorname{argmin}} \frac{||x-y||_2^2}{2} + \sigma F(y)$).
\begin{algorithm}[H]
\begin{algorithmic}
\STATE \textbf{Input: } Choose $\nu,\sigma > 0;\theta \in \{0,1\}; u^0 = \tilde u^0 \in X; y^0=0;$\\[.1in]
\STATE $k=0$\\[.1in]
\WHILE{Convergence criteria not met}
\STATE{$y^{k+1} = \operatorname{prox}_{\sigma F*}(y^k+\sigma K\tilde u^k)$ }\\[.1in]
\STATE {$u^{k+1} = \operatorname{prox}_{\nu G}(u^k-\nu K^* y^{k+1})$}\\[.1in]
\STATE {$\tilde u^{k+1} =u^{k+1} + \theta(u^{k+1} - u^k)$}\\[.1in]
\STATE $k=k+1$\\[.1in]
\ENDWHILE\\[.1in]
\STATE \textbf{Output: }Numerical solution $u$ to the saddle point problem.
\end{algorithmic}
\caption{General Primal-Dual algorithm.}
\label{alg:pd}
\end{algorithm}

Applied to the deconvolution problem given by \eqref{eq:deconvolution}, we have $K = D$, $F(.) = ||.||_1$ and $G(.) = \frac{\lambda}{2}||A\ast .-f||_2^2$. 

It is well known that the dual of $||.||_1$ is $||.||_\infty$ and furthermore, $\operatorname{prox}_{\sigma||.||_\infty}(x)$ corresponds to a projection of $x$ on to the unit $L^2$ ball \cite{L2projection}. Thus, denoting $(i,j)$ a pixel location,
\begin{equation}
\operatorname{prox}_{\sigma F*}(y)_{i,j} = \frac{y_{i,j}}{\max(1,|y_{i,j}|)},
\end{equation}
and
\begin{equation}
\operatorname{prox}_{\sigma F*}(y^k + \sigma D \tilde u^k)_{i,j} = \frac{y^k_{i,j} + \sigma (D \tilde u^k)_{i,j}}{\max(1,|y^k_{i,j} + \sigma D \tilde u^k)_{i,j}|)}.
\end{equation}

On the other-hand, $\operatorname{prox}_{\nu G}(z)$ is the minimization of two squared $L^2$ norms. Its solution has a closed-form formulation obtained by (we will denote $\mathbf{A}u=A\ast u$ to emphasize the operator point of view and where $\mathbf{A}^*$ denotes the adjoint of $\mathbf{A}$) 
\begin{equation} 
\operatorname{prox}_{\nu G}(z) = \underset{v}{\operatorname{argmin}} \frac{||z-v||_2^2}{2} +\nu\frac{\lambda}{2}||\mathbf{A} v-f||_2^2
\end {equation}
$$\implies 0 = -(z-v) + \nu \lambda \mathbf{A}^*(\mathbf{A}v-f)$$
$$\implies (I + \nu \lambda \mathbf{A}^*\mathbf{A})v = \nu \lambda \mathbf{A}^* f + z.$$
This can be easily solved in the Fourier domain by
\begin{equation}
\operatorname{prox}_\nu G(z) = \Big( (\widehat I + \nu \lambda |\widehat A|^2)^{-1} (\nu \lambda \overline{\widehat A} \widehat f + \widehat z)\Big)^\lor.
\end{equation}
Furthermore, since we are working in the finite dimensional image space, $D^* = D^T$. We then make the assumption that $D^T \approx D^{-1}$, and we numerically implement this using the inverse wavelet transform. Finally, the wavelet based deconvolution algorithm is given by Algorithm~\ref{alg:wavdeconv}.
\begin{algorithm}[H]
\begin{algorithmic}
\STATE \textbf{Input: } Choose $\nu,\sigma > 0;\theta \in \{0,1\}; x^0 = \tilde x^0 = f; y = Df;$\\[.1in]
\STATE $k=0$\\[.1in]
\WHILE{Convergence criteria not met}
\STATE{$y^{k+1}_{i,j} =  \frac{y^k_{i,j} + \sigma (D \tilde u^k)_{i,j}}{\max(1,|y^k_{i,j} + \sigma D \tilde u^k)_{i,j}|)}$}\\[.1in]
\STATE {$u^{k+1} = \Big( (\widehat I + \nu \lambda |\widehat A|^2)^{-1} (\nu \lambda \overline{\widehat A} \widehat f + \widehat u^k - \nu \widehat{D^{-1}y^{k+1}} )\Big)^\lor$}\\[.1in]
\STATE {$\tilde u^{k+1} =u^{k+1} + \theta(u^{k+1} - u^k)$}\\[.1in]
\STATE $k=k+1$\\[.1in]
\ENDWHILE\\[.1in]
\STATE \textbf{Output: }Numerical solution $u$ to the saddle point problem \eqref{eq:deconvolution}.
\end{algorithmic}
\caption{Wavelet based deconvolution.}
\label{alg:wavdeconv}
\end{algorithm}

\subsubsection{Extension to empirical wavelets}
\label{ch4:deconvewwt}
Here, we naturally aim at promoting sparsity in the empirical wavelets domain, i.e. solve the deconvolution problem \eqref{eq:deconvolution} by replacing the traditional wavelet operator, $D$ with an empirical wavelet operator, denoted $D_\epsilon$, as was done by Alvarado in \cite{zariluz}. Two strategies are possible: 1) we build the set of empirical wavelets based on the input image $f$ at the beginning, and this set is kept fixed while solving \eqref{eq:deconvolution}; 2) the set of wavelet filters is initialized based on $f$ at the beginning of the algorithm but then it is updated at each iteration. In order to distinguish both options, we will denote $D_\epsilon^0$ the empirical wavelet transform that keeps its set of filters fixed in option 1). It is worth noticing that in option 2), the operator $D_\epsilon$ is non-linear. Thus, we lose the theoretical guarantee of convergence of the primal-dual algorithm, however we never observed any convergence issues in our experiments. The algorithms corresponding to each options are given by Algorithms~\ref{alg:decew1} and \ref{alg:decew2}, respectively. 

\begin{algorithm}[H]
\begin{algorithmic}
\STATE \textbf{Input: } Choose $\nu,\sigma > 0;\theta \in \{0,1\}; u^0 = \tilde u^0 = f; $\\[.1in]
\STATE Build the set of empirical wavelet filters, $\{\widehat{\varphi_n}\}$ based on the magnitude spectrum of $f$. This defines the operator $D_\epsilon^0$\\[.1in]
\STATE Initialize $y^0=D_\epsilon^0 f,k=0$\\[.1in]
\WHILE{Convergence criteria not met}
\STATE{$y^{k+1}_{i,j} =  \frac{y^k_{i,j} + \sigma (D_\epsilon^0 \tilde u^k)_{i,j}}{\max(1,|y^k_{i,j} + \sigma D_\epsilon^0 \tilde u^k)_{i,j}|)}$}\\[.1in]
\STATE {$u^{k+1} = \Big( (\widehat I + \nu \lambda |\widehat A|^2)^{-1} (\nu \lambda \overline{\widehat A} \widehat f + \widehat u^k - \nu \widehat{(D_\epsilon^0)^{-1}y^{k+1}} )\Big)^\lor$}\\[.1in]
\STATE {$\tilde u^{k+1} =u^{k+1} + \theta(u^{k+1} - u^k)$}\\[.1in]
\STATE $k=k+1$\\[.1in]
\ENDWHILE\\[.1in]
\STATE \textbf{Output: }Numerical solution $u$ of $\min_u ||D_\epsilon^0 u||_1 + \frac{\lambda}{2}||Au-f||_2^2$.
\end{algorithmic}
\caption{Empirical wavelet based deconvolution with a kept fixed set of filters.}
\label{alg:decew1}
\end{algorithm}

\begin{algorithm}[H]
\begin{algorithmic}
\STATE \textbf{Input: } Choose $\nu,\sigma > 0;\theta \in \{0,1\}; u^0 = \tilde u^0 = y = f;$\\[.1in]
\STATE Initialize the set of empirical wavelet filters, $\{\widehat{\varphi_n^k}\}$ based on the magnitude spectrum of $f$.\\[.1in]
\STATE $k=0$\\[.1in]
\WHILE{Convergence criteria not met}
\STATE{$y^{k+1}_{i,j} =  (D_\epsilon^k)^{-1}\frac{D_\epsilon^k y^k_{i,j} + \sigma (D_\epsilon^k \tilde u^k)_{i,j}}{\max(1,|D_\epsilon^k y^k_{i,j} + \sigma D_\epsilon^k \tilde u^k)_{i,j}|)}$}\\[.1in]
\STATE{$u^{k+1} = \Big( (\widehat I + \nu \lambda |\widehat A|^2)^{-1} (\nu \lambda \overline{\widehat A} \widehat f + \widehat u^k - \nu \widehat y^{k+1})\Big)^\lor$}\\[.1in]
\STATE {$\tilde u^{k+1} =u^{k+1} + \theta(u^{k+1} - u^k)$} \\[.1in]
\STATE Update the set of empirical filter $\{\widehat{\varphi_n^k}\}$ based on the magnitude spectrum of $u^{k+1}$\\[.1in]
\STATE $k=k+1$\\[.1in]
\ENDWHILE \\[.1in]
\STATE \textbf{Output: }Numerical solution $u$.
\end{algorithmic}
\caption{Empirical wavelet based deconvolution with sets of wavelet filters that are updated at each iteration.}
\label{alg:decew2}
\end{algorithm}

We want to emphasize that updating the set of empirical wavelet filters $\{\widehat{\varphi_n^k}\}$ at each iteration $k$ automatically updates the operator $D_\epsilon^k$.

\subsubsection{Results}
In this section, we compare the deconvolution results obtained by using framelets \cite{framelets}, empirical curvelet option 1 and the empirical watershed wavelets. We use the Outex dataset once again, and simulate a Gaussian blur with variance $\sigma_B^2 = 1,2,3$. To assess the effectiveness of the different restorations, we compute the structural similarity index measure (SSIM) 
\cite{Wang2004} between the pristine image $f$ and the deblurred image $u$. We remind the reader that the closer the SSIM to 1, the better the deconvolution. The deconvolution algorithms are run on 50 images, and we present the average SSIM over these 50 images. 

We begin by investigating the influence of the parameters of the EWWT in a similar way as we did in Section~\ref{ch4:ewwtparams}. We run these tests in the case of the fixed filter bank approach described in Section~\ref{ch4:deconvewwt}. We first consider the influence of the scale-space step-size $s_0$ while keeping $\tau = 0.1$. The results are presented in Table~\ref{tab:deconvstepsize}.
 
 \begin{table}[H]
\caption{Influence of the scale-space step-size $s_0$ for deconvolution for different level of blur ($\sigma_B^2$). The best SSIM values are given in bold}.\label{tab:deconvstepsize}
\centering
\begin{tabular}{cccccccccccccccccccccc}
\toprule
$\boldsymbol{s_0}$ &  \textbf{0.05} &  \textbf{0.1} &  \textbf{0.2} &  \textbf{0.3} &  \textbf{0.4} &  \textbf{0.5}  &  \textbf{0.6} &  \textbf{0.7} &  \textbf{0.8} &  \textbf{0.9} &  \textbf{1.0}\\
\midrule
$\sigma_B^2 = 1$ & \textbf{0.9931} & \textbf{0.9931} & 0.9930 & 0.9930 & 0.9930 & 0.9930 & 0.9930 & 0.9927 & 0.9927 & 0.9926 & 0.9925 \\

$\sigma_B^2 = 2$ & \textbf{0.9305} & 0.9270 & 0.9120 & 0.9072 & 0.9000 & 0.8960 & 0.8918 & 0.8882 & 0.8877 & 0.8866 & 0.8863\\

$ \sigma_B^2 = 3$ & 0.7554 & \textbf{0.7556} & 0.7551 & 0.7550 & 0.7549 & 0.7549 & 0.7533 & 0.7522 & 0.7487 & 0.7470 & 0.7400\\
\bottomrule
\end{tabular}
\end{table}

We observe that a scale-space step-size value of $s_0 = 0.05$ or $s_0 = 0.1$ is optimal. While we see small differences between average SSIMs, we notice that smaller step-sizes give higher average results, confirming our observation from Section~\ref{ch4:ewwtparams} that a larger number of modes is advantageous. It is worth noticing that the value $s_0 = 0.1$ also led to the best texture segmentation, thus we favor the value $s_0=0.1$ over the value $s_0 = 0.05$ in the following experiments.

Next, we compare the influence of the transition width $\tau$. From Table~\ref{tab:deconvtau}, we see very small differences between the structural similarity index measures, but overall a small transition width of $\tau = 0.05$ works best for less blur, while a larger value of $\tau = 0.5$ works best for larger blurs. Going forward, we use $\tau = 0.05$ for $\sigma_B^2 = 1,2$ and $\tau = 0.5$ for $\sigma_B^2 = 3$.

\begin{table}[H]
\caption{Influence of transition width $\tau$ for deconvolution for different level of blur ($\sigma_B^2$). The best SSIM values are given in bold.}\label{tab:deconvtau}
\centering
\begin{tabular}{ccccccccc}
\toprule
$\tau$ & \textbf{0.05} &  \textbf{.1} &  \textbf{.2} & \textbf{ .3} &  \textbf{.5} \\
\midrule
$\sigma_B^2=1$ & \bf{0.9931} & \bf{0.9931} & 0.9930 & 0.9930 & .9930 \\
$\sigma_B^2=2$ & \bf{0.9276} & 0.9270 & 0.9270 & 0.9268 & 0.9266 \\
$\sigma_B^2=3$ & 0.75586 & 0.7556 & 0.7559 & 0.7560 & \bf{0.7564} \\
\bottomrule
\end{tabular}
\end{table}

We now compare the effectiveness of empirical watershed wavelet deconvolution against framelets and empirical curvelet option 1 on the Outex dataset. Framelets have been chosen as the classic wavelets since it has been shown to outperform other traditional wavelet decompositions \cite{cai1}, \cite{cai2}, \cite{cai3}.

\begin{figure}[!t]
\centering
\subfloat[Original]{\includegraphics[width=.22\textwidth]{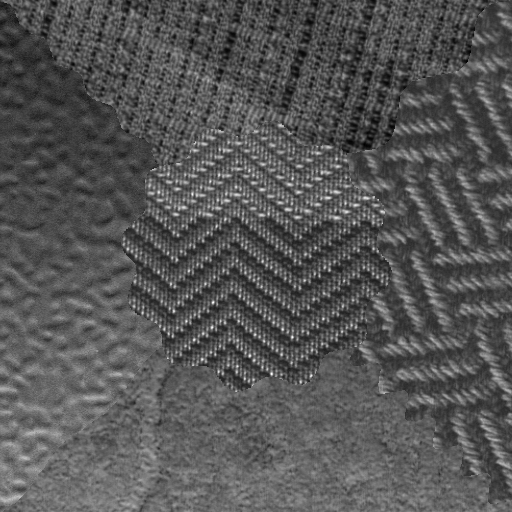}
\includegraphics[width=.22\textwidth]{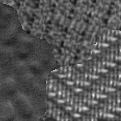}\label{77orig}} \hspace{.02in}
\subfloat[Blurred $\sigma_B^2 = 3$]{\includegraphics[width=.22\textwidth]{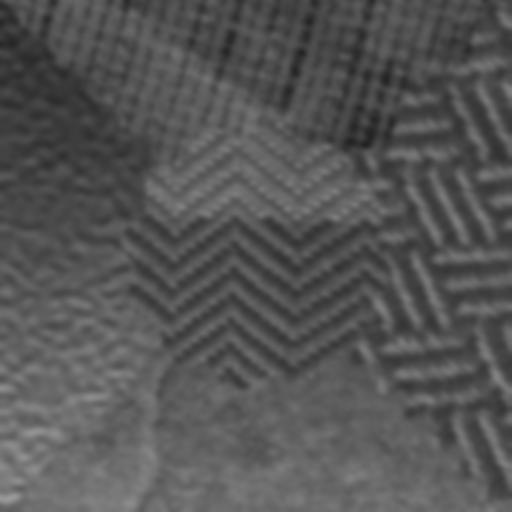}
\includegraphics[width=.22\textwidth]{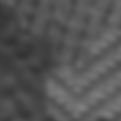}\label{77blurred}}\\ 
\subfloat[EWT-C1 Fixed]{\includegraphics[width=.22\textwidth]{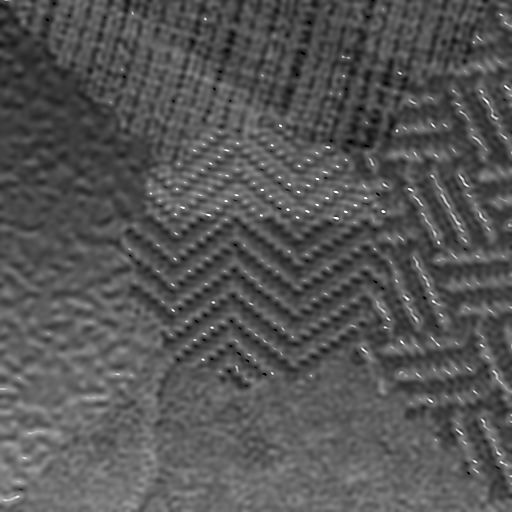}  
\includegraphics[width=.22\textwidth]{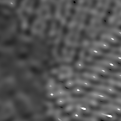}\label{77ewtc1}}  \hspace{.02in}
\subfloat[EWT-C1 Adaptive]{\includegraphics[width=.22\textwidth]{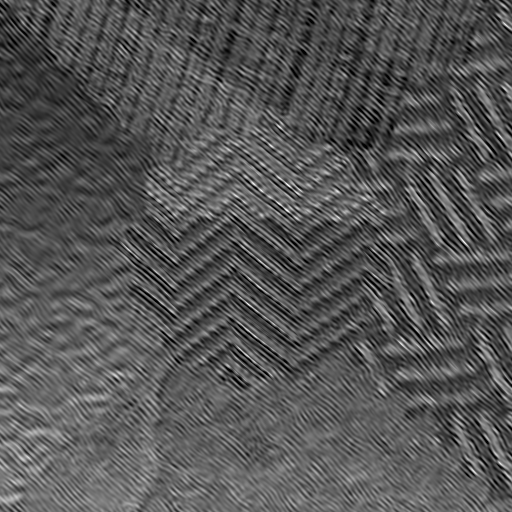}  
\includegraphics[width=.22\textwidth]{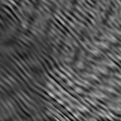}\label{77ewtcn}} \\
\subfloat[Framelet]{\includegraphics[width=.22\textwidth]{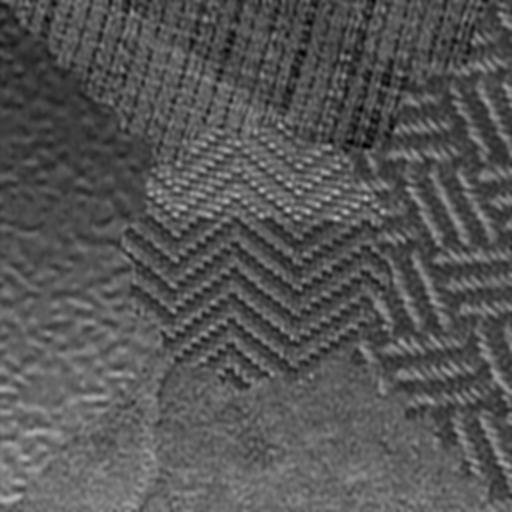} 
\includegraphics[width=.22\textwidth]{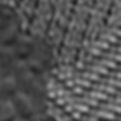}\label{77framelet}} \hspace{.02in}
\subfloat[EWWT Fixed]{\includegraphics[width=.22\textwidth]{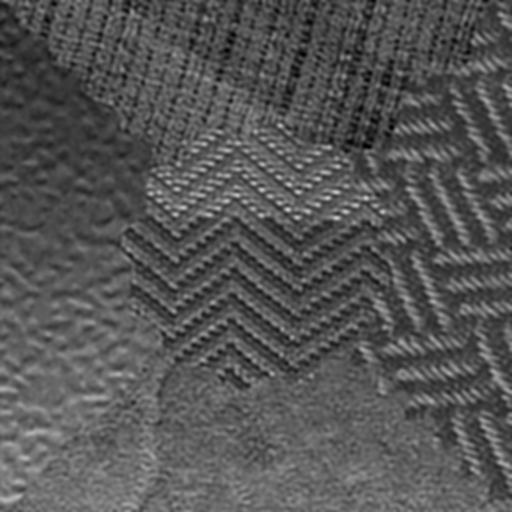} 
\includegraphics[width=.22\textwidth]{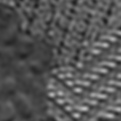}\label{77ewwt1}} \\
\subfloat[EWWT Adaptive]{\includegraphics[width=.22\textwidth]{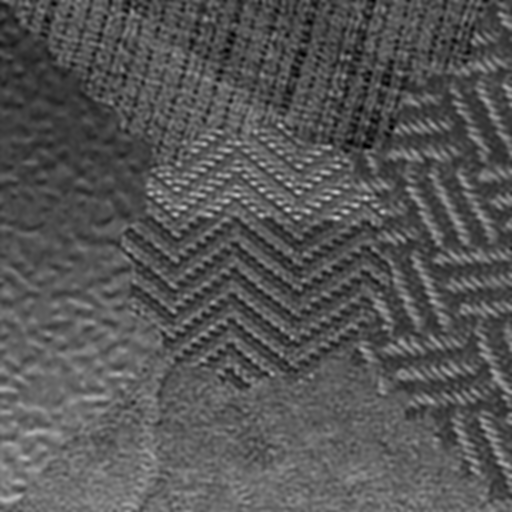}  
\includegraphics[width=.22\textwidth]{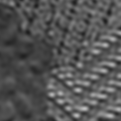}\label{77ewwtn}} 
\caption{Example of the deconvolution performed on an image with a blur variance of $\sigma_B^2 = 3$. Note the artifacts present in the empirical curvelet approaches.\label{fig:deconvcompare1}}
\end{figure}

\begin{figure}[!t]
\centering
\subfloat[Original]{\includegraphics[width=.22\textwidth]{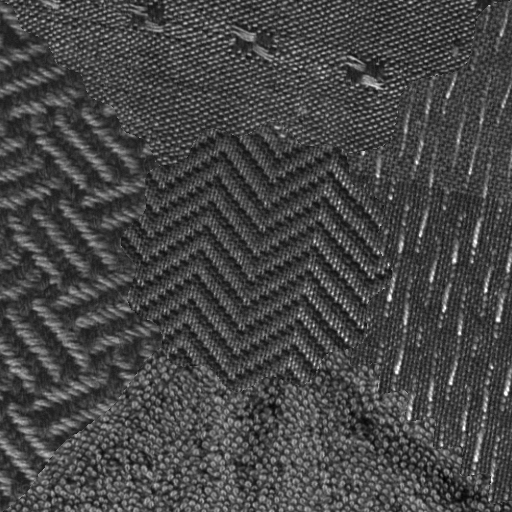}
\includegraphics[width=.22\textwidth]{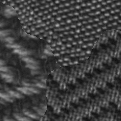}\label{81orig}} \hspace{.02in}
\subfloat[Blurred $\sigma_B^2 = 2$]{\includegraphics[width=.22\textwidth]{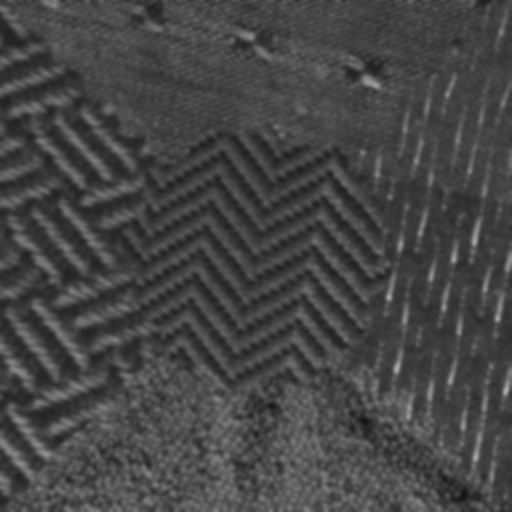}
\includegraphics[width=.22\textwidth]{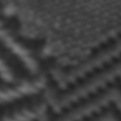}\label{81blurred}}\\ 
\subfloat[EWT-C1 Fixed]{\includegraphics[width=.22\textwidth]{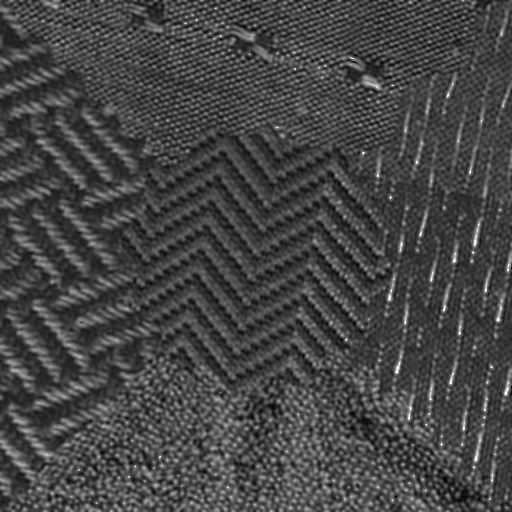}  
\includegraphics[width=.22\textwidth]{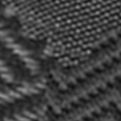}\label{81ewtc1}}  \hspace{.02in}
\subfloat[EWT-C1 Adaptive]{\includegraphics[width=.22\textwidth]{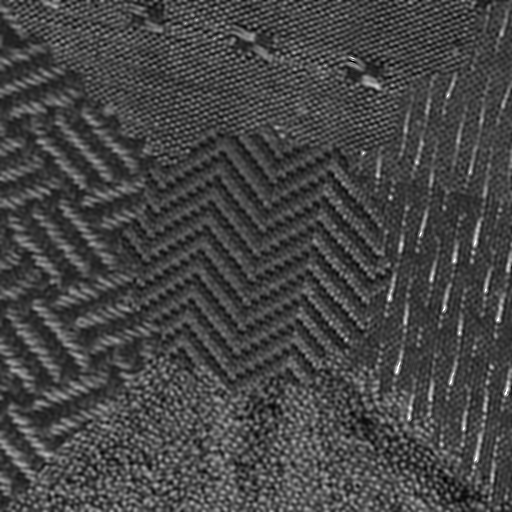}  
\includegraphics[width=.22\textwidth]{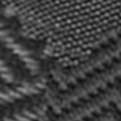}\label{81ewtcn}} \\
\subfloat[Framelet]{\includegraphics[width=.22\textwidth]{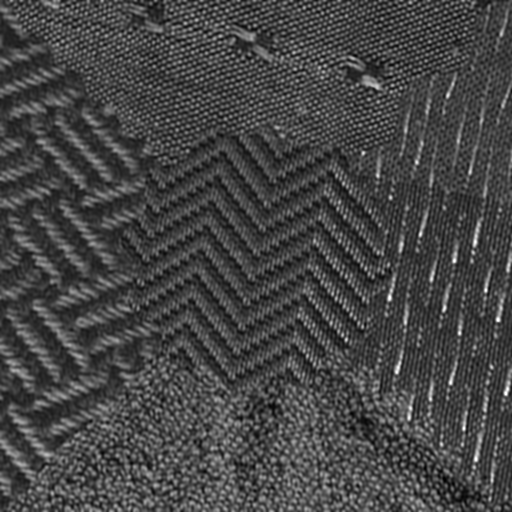} 
\includegraphics[width=.22\textwidth]{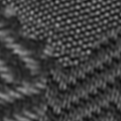}\label{81framelet}} \hspace{.02in}
\subfloat[EWWT Fixed]{\includegraphics[width=.22\textwidth]{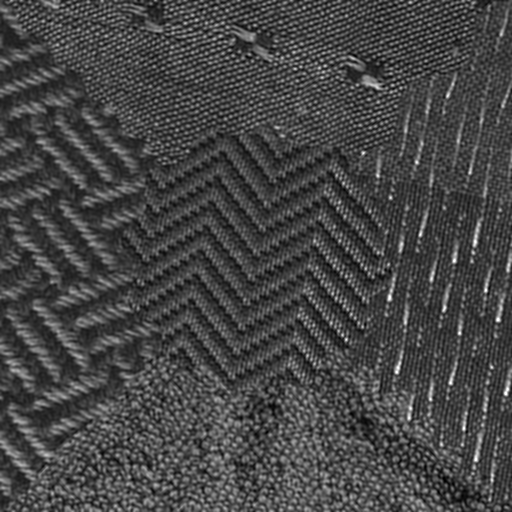} 
\includegraphics[width=.22\textwidth]{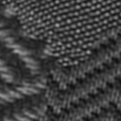}\label{81ewwt1}} \\
\subfloat[EWWT Adaptive]{\includegraphics[width=.22\textwidth]{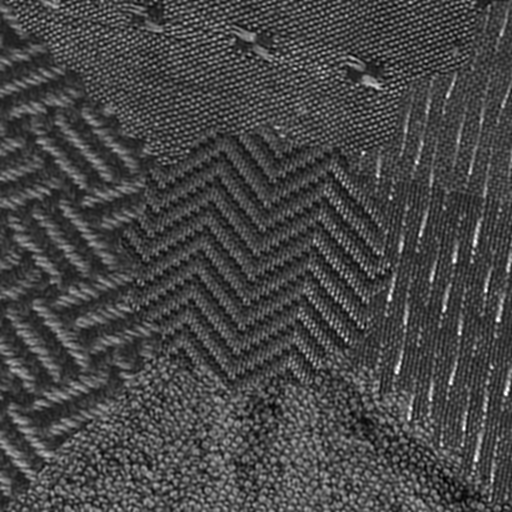}  
\includegraphics[width=.22\textwidth]{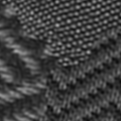}\label{81ewwtn}} 
\caption{Example of the deconvolution performed on an image with a blur variance of $\sigma_B^2 = 2$.\label{fig:deconvcompare2}}
\end{figure}

For both empirical wavelet approaches, we look into both the fixed filter bank approach  and the adaptive filter bank method described in Section~\ref{ch4:deconvewwt}. The primal-dual algorithm has its own parameters set to $\lambda = 5000$, $\sigma = 0.1$, $\nu = 0.1$, and $\theta = 1$ \cite{zariluz}. The results (again averaged on 50 images) are shown for all $\sigma_B^2$ values in Table~\ref{tab:deconvcomparison}, while visual deblurring examples are shown in Figure~\ref{fig:deconvcompare1} and \ref{fig:deconvcompare2}. 

\begin{table}[H]
\caption{Final comparisons for deconvolution for different level of blur ($\sigma_B^2$). The best SSIM values are given in bold.}\label{tab:deconvcomparison}
\centering
\begin{tabular}{ccccccccccc}
\toprule
\textbf{Method} & \textbf{Framelets} & \multicolumn{2}{c}{\textbf{EWT-C1}} & \multicolumn{2}{c}{\textbf{EWWT}} \\
&  &  \textbf{fixed} & \textbf{adaptive} & \textbf{fixed} &\textbf{adaptive} \\
\midrule
$\sigma_B^2=1$ & 0.9929 & 0.9922 & 0.9701 & \textbf{0.9931} & \textbf{0.9931} \\
$\sigma_B^2=2$ & 0.9172 & 0.8698 & 0.8449 &  0.9276 & \textbf{0.9337} \\
$\sigma_B^2=3$ & 0.6922 & 0.6362 & 0.5973 &  0.7564 & \textbf{0.7575} \\
\bottomrule
\end{tabular}
\end{table}

We see from these results that among all tested methods, the empirical watershed wavelet performs the best. Furthermore, between the fixed filter bank and the adaptive filter bank approaches, the adaptive method slightly outperforms the fixed method. We also note that in some cases, the empirical curvelets led to artifacting, which can be seen in Figure~\ref{fig:deconvcompare1}, and this phenomenon was not observed when using framelets and empirical watershed wavelets. In Table~\ref{tab:deconviters}, we compare the average number of iterations needed to reach the desired convergence accuracy (set to $2\times 10^{-4}$) between the framelet and the empirical watershed wavelet methods. We see that the empirical watershed wavelets converge in comparable or less iterations compared to other methods. This seems to be another advantage of the adaptability, i.e. adaptive representations need less iterations to converge to optimal representations of the images.

\begin{table}[H]
\caption{Average number of iteration until convergence.\label{tab:deconviters}}
\centering
\begin{tabular}{ccccccccccc}
\toprule
\textbf{Method} & \textbf{Framelets} &  \multicolumn{2}{c}{\textbf{EWWT}} \\
&  &  \textbf{fixed} & \textbf{adaptive} \\
\midrule
$\sigma_B^2=1$ & \textbf{7.92} &   8.1 & 8.1 \\
$\sigma_B^2=2$ & 31.18 &  32.7 & \textbf{26.62}\\
$\sigma_B^2=3$ & 45.7 &  34.2 & \textbf{33.2} \\
\bottomrule
\end{tabular}
\end{table}

\section{Conclusion}
\label{sec:5}
In conclusion, we proposed a new construction of the 2D empirical wavelet transform based on an arbitrary partitioning of the Fourier spectrum. This construction includes a trigonometric transition region similar to the 1D EWT, and we showed that the set of empirical wavelets forms a frame. We also proposed an adaptive way to obtain such a partitioning, which is broken up into two steps. The first step involves a mode detection step and uses scale-space representations to define persistent maxima. The second step is a boundary detection step that uses the watershed transform. This defines a boundary that isolates the detected maxima by a path of lowest separation. Using these in conjunction with the 2D EWT of arbitrary shape, we defined the empirical watershed wavelet transform (EWWT). 

We then showed that the EWWT performs well when applied to multiple applications. Both as a feature extractor for unsupervised texture segmentation and as a sparsity promoter in deconvolution, we saw improved results over previous wavelet and empirical wavelet transforms. On top of this, we explored the robustness of the empirical watershed wavelet transform parameters on each of these applications and found trends in how these parameters act.

Future work is possible in multiple directions. Firstly, in the construction of the 2D EWT, one could consider an arbitrary shape construction that forms a tight frame, rather than a frame. Future work can also be considered in the extension from 2D to 3D or N-D versions of the EWT of arbitrary shape. 

In terms of the mode detection we proposed, many other options exist, such as threholding maxima, using maxima filters, entropy-based methods like in wavelet packets, or PCA based methods for dimension reduction. There is plenty of work to be done in finding which method works best for various applications. The same can be said for the boundary detection, especially if certain degrees of boundary smoothness are desired. While we saw that a certain value of the scale-space step-size worked best for both applications explored, more work is warranted to see if this extends to further applications, as well as further datasets. 

\section{Acknowledgement}
This work was partially supported by the Air Force Office of Scientific Research under the grant number FA9550-15-1-0065.

\reftitle{References}




\end{document}